\documentclass[12pt]{article}
\usepackage{amsfonts}
\usepackage{amssymb}
\usepackage{mathtools}
\usepackage{verbatim}
\usepackage{amsthm}
\theoremstyle{plain}
\usepackage[english, russian]{babel}

\newtheorem*{theorem*}{Theorem}

\newtheorem{definition*}{Definition}

\begin{document}

\begin{center}
{\Large \bf On optimal packing of Minkowski spheres
}
\end{center}


\begin{center}
{\bf Nikolaj M. Glazunov } \end{center}

\begin{center}

{\rm Glushkov Institute of Cybernetics NASU, Kiev, } \\

{\rm  Email:} {\it glanm@yahoo.com }
\end{center} 

{\bf 2020 Mathematics Subject Classification:}  11H06, 11H31, 11Jxx, 52A15, 52C05, 90C25 \\ 

\bigskip

{\bf Keywords:} Minkowski sphere, Minkowski-Cohn surface,   critical lattice, optimal lattice packing, packing density.\\

\textbf{Annotation}\\

Elements of the theory of non-differentiable ravine-type functions and the theory of a corresponding class of non-differentiable ravine-type manifolds—specifically, ravine-type manifolds with corners—have been developed. A continuum class of non-differentiable, essentially ravine-type functions has been constructed. A non-differentiable ravine-type function has been formulated, the corresponding manifold of which is a non-differentiable Minkowski-Cohn surface. \\

Based on the non-differentiable ravine-type Minkowski-Cohn surface, the author has developed a method for the optimal packing of a continuum family of one-dimensional Minkowski spheres in the real plane. By applying this method to the packing of one-dimensional Minkowski spheres in the plane, an optimal lattice packing for a continuum family of such spheres has been obtained for the first time. Expressions for these optimal packing lattices have been derived, and their validity has been proven. Formulas for the packing density corresponding to these lattices have also been established.\\

\textbf{Introduction}\\

The study and solution of optimization problems are conducted within the framework of various theories and methods. A primary challenge in developing efficient (sub)gradient methods for convex optimization stems from the ill-conditioning of significant classes of applied problems: these problems involve functions with highly elongated level sets (often referred to as "ravine-like" or "valley-like" functions). This issue is particularly relevant to non-smooth optimization problems—characterized by "pronouncedly ravine-like" functions—and problems involving singular features, such as the degeneracy of the functional's Hessian at the optimal point, linear dependence of subgradients, a degenerate structure of the solution set, or non-uniqueness of the solution. The development of numerical optimization methods is largely driven by the need to overcome these challenges. While this enables the derivation of approximate optimal solutions, such methods may prove insufficient for obtaining exact optimal solutions.

The application of "ravine functions" (or "valley functions") to optimization problems was proposed by I. M. Gelfand and M. L. Tsetlin \cite{src1} in 1961 and further developed by N. Z. Shor and other researchers \cite{src2,src3,src4}.

The author learned about these methods in the 1970s while undertaking research internships at the Computing Center of the USSR Academy of Sciences, Moscow State University, and the Steklov Mathematical Institute of the USSR Academy of Sciences, where he participated in scientific seminars. These internships were arranged at the author's request by Academician V. M. Glushkov, Director of the Institute of Cybernetics. The first internship, lasting approximately one and a half years, took place in 1970–1971 during the author's postgraduate studies at the Institute of Cybernetics; the second internship, also lasting about one and a half years, occurred in 1974–1975, following the author's successful defense of his Candidate of Sciences dissertation.

One of the characteristics of ravine functions is that their level curves are elongated. Let us consider a class of real-valued ravine functions that are either convex or concave. Concave functions with highly elongated level surfaces (and corresponding level curves) could be described as mountainous or as height functions; however, we shall also refer to such functions as ravine functions. Yaruzhni (ridge-like/valley-like) functions arise in problems within number theory, approximation theory, cybernetics, and discrete analysis, as well as in issues related to the control processes of physiological mechanisms. Optimization methods proposed by Gelfand and Tsetlin--along with subsequent methods developed by them and their followers--pertain to the numerical optimization of a function at a specific point. In this context, we should mention the space-transformation methods developed within the school of N. Z. Shor by his students (now researchers at the Glushkov Institute of Cybernetics of the National Academy of Sciences of Ukraine), including M. G. Zhurbenko and others \cite{src2,src3,src4}. However, a range of problems in number theory, approximation theory, cybernetics, discrete analysis, and applied 

This text is devoted to the investigation and potential solution methods for such a problem, using the example of a function whose range is a Minkowski-Cone surface. Although a hint regarding the construction of such a surface appeared in Minkowski’s well-known monograph \cite{src5}, and the surface itself—in terms of a $(\tau, \sigma)$ parametrization—was constructed by Cohn \cite{src6} in 1950, the study of this surface and the problems associated with it continues to attract the attention of researchers. In particular, to the best of the author’s knowledge, the curve of maxima on this surface remains unknown. Recently (2025–2026), researchers from Austria and Australia contacted the author regarding results and references to the author’s work in this field, and the author provided the relevant references.

The content of this text is as follows. We begin by presenting elements of the geometry and analysis of "ravine" functions. We define the function that characterizes the Minkowski–Cohn surface and develop elements of optimization theory for the first parametrization of this surface—known as the Cohn parametrization. When determining the curve of minima on this surface, it is convenient to employ different parametrizations depending on the specific subregion of the Minkowski–Cohn surface under consideration. For established results, see \cite{src7, src8, src9, src10} and the references cited therein.

\bigskip

\textbf{1.1. On Stratified Functions}\\

\textbf{1.1.1. On Non-differentiable Functions and Stratified-type Manifolds}\\

We shall consider non-differentiable functions of a stratified type and the corresponding class of stratified non-differentiable manifolds—namely, stratified manifolds with corners. Since such functions and manifolds constitute an extension of the class of smooth (differentiable) functions and manifolds, we begin by recalling the properties of smooth objects that we will need. We shall consider functions defined on closed and open subsets of the Euclidean space \({\mathbb{R}}^{n}\).\\

\emph{\textbf{Example 1}}. The set \([a, b]\) and the set \((a, b)\) are, respectively, a closed interval and an open interval in \(\mathbb{R}\).

The set \([a_1, b_1] \times \dots \times [a_n, b_n] \subset \mathbb{R}^n\) is called a closed rectangle in \(\mathbb{R}^n\), and the set \((a_1, b_1) \times \dots \times (a_n, b_n) \subset \mathbb{R}^n\) is called an open rectangle.

A set \(U \subset \mathbb{R}^n\) is called open if, for every \(x \in U\), there exists an open rectangle \(R\) such that \(x \in R \subset U\). A subset \(V\) of \(\mathbb{R}^n\) is closed if \(\mathbb{R}^n \setminus V\) is open.

\textbf{Definition 1}. Let \(E\) and \(F\) be topological spaces. Let \(f\) be a one-to-one mapping from \(E\) onto \(F\). If both \(f\) and its inverse mapping are continuous, then \(f\) is called a homeomorphism, and the spaces \(E\) and \(F\) are said to be homeomorphic.

\textbf{1.1.2 On Smooth Functions and Smooth Manifolds}\\

Let us recall some concepts and results. In what follows, we refer to \cite{src11,src12,src13,src14}.\\

\emph{\textbf{Definition 2}. Let \(U\) be an open \(n\)-dimensional set in the Euclidean space \({\mathbb{R}}^{n}\), and let \(f\) be a real-valued function defined on \(U\). The function \(f\) is called smooth if, at every point in \(U\), it possesses continuous partial derivatives of all orders with respect to the coordinates in \({\mathbb{R}}^{n}\). In this case, \(f\) is called a smooth mapping of the set \(E\) into \({\mathbb{R}}^{n}\).}\\

\emph{\textbf{Definition 3}}. Let \(E\) be a set in \({\mathbb{R}}^{n}\), and let \(f_{1}, f_{2}, \dots, f_{n}\) be \(n\) smooth functions on \({\mathbb{R}}^{n}\). The mapping $f(x)$ of a set $E \subset \mathbb{R}^n$ into $\mathbb{R}^n$ assigns the point $(f(x_1), \dots, f(x_n))$ as the image of the point $(x_1, \dots, x_n)$.\\

An open $n$-dimensional ball in the Euclidean space $\mathbb{R}^n$ is the set\\

\begin{equation}\label{eq:open-ball}B^{n}\  = \ \{(x_{1},\ .\ .\ .\ ,\ x_{n})|x_{1}^{2}\  + \  \cdot \  \cdot \  \cdot \ x_{n}^{2}\  < \ 1\}.\end{equation}\\

A set $W$ in the Euclidean space $\mathbb{R}^n$ is called an open cell if it is homeomorphic to the open ball \eqref{eq:open-ball}.\\

\emph{\textbf{Definition 4}. An n-dimensional smooth manifold M --- is a Hausdorff topological space with a covering consisting of a countable number of open sets \(U_{1}, U_{2}, \dots\) that satisfy the following conditions:}\\

\emph{1) For each \(U_{i}\), there exists a homeomorphism \(\phi_{i}: U_{i} \to W_{i}\), where \(W_{i}\) is an open set in Euclidean space.}

\emph{2) If \(U_{i} \cap U_{j} \neq \varnothing\), then the homeomorphism (transition map) \(\phi_{ij} = \phi_{j} \circ \phi_{i}^{-1}\) from the set \(\phi_{i}(U_{i} \cap U_{j})\) to \(\phi_{j}(U_{i} \cap U_{j})\) is a smooth map.}

\emph{The homeomorphisms \(\phi_{i}: U_{i} \to W_{i}\) are called coordinate charts, and the collection \(\{U_{i}, \phi_{i}\}\) is called an atlas of M.}

\emph{The number n is called the dimension (dim M) of M (an n-manifold).}\\

\textbf{1.1.3 Manifolds with Boundary}\\

Let

\begin{equation}\label{eq:upper-halfspace}H^{m}\  = \ \{ x\  = \ (x_{1},\ .\ .\ .\ ,\ x_{m})\  \in \ {\mathbb{R}}^{m}|x_{m}\  \geq \ 0\},\ \ \end{equation}\\

and

\begin{equation}\label{eq:halfspace-boundary}\partial H^{n}\  = \ \{ x\  = \ (x_{1},\ .\ .\ .\ ,\ x_{m})\  \in \ {\mathbb{R}}^{m}|x_{m}\  = \ 0\}\ \end{equation}

be the \emph{m}-dimensional upper half-space and its boundary.\\

\textbf{\emph{Definition} 5}. A smooth $m$-manifold with boundary consists of a Hausdorff topological space $M$, an open covering $\{U_i\}$ of $M$, and a collection of homeomorphisms

\[\phi_i : U_i \rightarrow \mathfrak{U}_i\]

onto open subsets $\mathfrak{U}_i \subset H^m$—one for each $i$—such that for every pair $i, j$, the transition map

\[\phi_{ij} = \phi_j \phi_i^{-1} : \phi_i(U_i \cap U_j) \rightarrow \phi_j(U_i \cap U_j)\]

is a diffeomorphism.

A point $p \in M$ is called a boundary point if it satisfies the following condition: if $p \in U_i$, then $\phi_i(p) \in \partial H^m$.

The set

\begin{equation}\label{eq:boundary-of-manifold}\partial M = \{ p \in M \mid \phi_i(p) \in \partial H^m \text{ for all } i \text{ such that } p \in U_i \}\end{equation}

of all boundary points is called the boundary of $M$.\\

\textbf{1.1.4 Manifolds with corners}\\ 

Convex combination of points \(x_{1}, \dots, x_{n}\) in \({\mathbb{R}}^{n}\) is a linear combination

\begin{center}
{\(\displaystyle 
 \lambda_{1}x_{1} + \dots + \lambda_{n}x_{n}, \quad \text{where } \lambda_{1} + \dots + \lambda_{n} = 1 \quad \text{and } \lambda_{1}, \dots, \lambda_{n} \geq 0.\)}\end{center}

\textbf{Theorem}. A subset \(C\) of \({\mathbb{R}}^{n}\) is convex if and only if any convex combination of points from \(C\) also lies in \(C\).

Let

\begin{equation}
\displaystyle H^{n,k} = [0, \infty)^{k} \times {\mathbb{R}}^{n - k} = \{ x \in {\mathbb{R}}^{n} \mid x_{i} \geq 0, 1 \leq i \leq k \},
\end{equation}

Using these sets \eqref{eq:corner-halfspace}, one can define manifolds with corners (analogously to the use of \eqref{eq:upper-halfspace} for manifolds with boundary).\\

\textbf{1.1.5 On a Class of Gully-Shaped Functions}\\

\emph{Shor's monograph \cite{src2} considers a class of convex functions $f$ for which, using the notation from \cite{src2}, the following inequality holds:}

\begin{equation}\label{eq:shor-ravine-inequality}(g_{f}\ (x),\ x\  - \ x^{\ast}(x))\  \geq \ cos\ \varphi||g_{f}\ (x)||\ ||x\  - \ x^{\ast}(x)||,\end{equation}

\emph{Here, $f$ is a convex function defined on $\mathbb{R}^n$, $g_f(x)$ is a subgradient of $f$, $0 \leq \varphi < \pi/2$, $x \in \mathbb{R}^n$, and $x^*(x)$ is the point in the set of minimizers of $f$ that is closest to $x$. If there exists no $\varphi < \pi/2$ such that \eqref{eq:shor-ravine-inequality} holds for all $x$ in a neighborhood of an optimal point of $f$, then $f$ is said to be essentially gully-shaped (or essentially gully-like).}

\emph{\textbf{Remark.} In Shor's monograph \cite{src2}, such functions are referred to as "essentially gully-shaped functions."} This name is actually more fitting for the functions we are investigating; however, we will use the term "valley-like" (*yaruzhno-podibnyi*), as it is more commonly employed.

An example of a function that is "essentially valley-like" is the one proposed by M. G. Zhurbenko \cite{src4}: the maximum of two positive binary quadratic forms. Below, we present an infinite set of essentially valley-like functions, constructed based on shifted one-dimensional Minkowski spheres.

\emph{\textbf{Example 2}.}

\begin{equation}\label{eq:zhurbenko-maximum}F_{p}\  = \ max\{ f_{p},\ g_{p}\},\ p\  > \ 1.\end{equation}

\begin{equation}\label{eq:zhurbenko-forms}f_{p}\  = \ |x\  + \ 1|^{p}\  + \ |y|^{p},\ g_{p}\  = \ |x\  - \ 1|^{p}\  + \ |y|^{p}.\end{equation}

\emph{\textbf{Remark 1.} We shall refer to the essentially valley-like functions in Example 2 as Zhurbenko's essentially valley-like functions (or simply Zhurbenko's valley-like functions).}

\emph{It can be conjectured that for $m$-dimensional Minkowski spheres (where $m > 1$), one can also construct a continuous family of essentially valley-like functions.}

\emph{A similar conjecture can be made regarding the set of Minkowski domains \cite{src25,src27}.}\\

\textbf{1.1.6 The Minkowski-Cohn Surface in the Cohn Parameterization.}\\

The "ravine" functions described below (as well as those described in subsequent sections of this work) are related to the study of certain problems in approximation theory, coding theory, and covering problems. In what follows, we employ notation and results from \cite{src5,src6,src7,src8,src9,src10}.\\

\emph{Let us define the concepts used. We consider balls of the form}

\[D_{p}\ :\ |x|^{p}\  + \ |y|^{p}\  \leq \ 1,\ p\  \geq \ 1,\]

\emph{and refer to such balls as Minkowski balls.} Continuing this, we consider the following classes of Minkowski balls and circles (one-dimensional spheres):\\

\emph{Watson balls:} $|x|^p + |y|^p \leq 1$ for $p_0 > p \geq 2$;\\

\emph{Davis balls:} $|x|^p + |y|^p \leq 1$ for $p_0 > p \geq 2$;\\

\emph{Mordell–Chebyshev balls:} $|x|^p + |y|^p \leq 1$ for $p_0 > p \geq 2$;\\

\emph{\textbf{Remark 1}. This classification corrects the one presented in \cite{src6} and more accurately reflects the historical contributions of the researchers.}\\

We present an analytical parametrization, following Cohn, of the admissible lattices for these domains, which include the critical lattices. It follows from \cite{src22} that such admissible lattices each possess 3 pairs of points on the boundary of the domain \(D_{p}\). Let \(x_{1} > x_{2} \geq 0\), \(y_{1} \geq 0\), and \(y_{2} \geq 0\) be real numbers such that the vectors \((x_{1}, y_{1})\) and \((-x_{2}, y_{2})\) form a basis for a \(D_{p}\)-admissible lattice that has three pairs of points on the boundary of the domain \(D_{p}\). Let us introduce real parameters \(\sigma \geq \tau \geq 0\) and set

\begin{center}
{\(\displaystyle x_{1}\  = \ \frac{1}{(1 + \tau^{p})^{1/p}},\ y_{1}\  = \ \frac{\tau}{(1 + \tau^{p})^{1/p}},\ x_{2}\  = \ \frac{1}{(1 + \sigma^{p})^{1/p}},\ y_{2}\  = \ \frac{\sigma}{(1 + \sigma^{p})^{1/p}}\ .\)}\end{center}

Using this notation and the concept of a determinant, we obtain the equation

\begin{equation}\label{eq:minkowski-cohn-surface}
\displaystyle \Delta(p,\sigma)=(\tau+\sigma)(1+\tau^{p})^{-1/p}(1+\sigma^{p})^{-1/p}
\end{equation}
which has a limiting point at \(p=1\), parameterizes lattices, and is defined within the region \(\gamma\): \(1 \leq p\), \(1 \leq \sigma \leq \sigma_{p} = (2^{p} - 1)^{1/p}\) in the \(\{p, \sigma\}\) plane, where \(\sigma\) is a real parameter; Here, \(\tau = \tau(p, \sigma)\) is a function uniquely determined by the condition

\[A^p + B^p = 1, \quad 0 \leq \tau \leq \tau_p,\]

where

\[A = A(p, \sigma) = (1 + \tau^p)^{-1/p} - (1 + \sigma^p)^{-1/p},\]

\[B = B(p, \sigma) = \tau(1 + \tau^p)^{-1/p} + \sigma(1 + \sigma^p)^{-1/p},\]

and \(\tau_p\) is determined by the equation

\[2(1 - \tau_p)^p = 1 + \tau_p^p, \quad 0 < \tau_p < 1.\]

We shall call this surface the Minkowski-Cohn surface. We have the following result:

\textbf{Proposition 1.} \emph{The Minkowski-Kohn surface parametrizes all admissible lattices of a domain that have three pairs of points on the boundary of this domain.}\\

\textbf{1.1.7. Elements of the analysis of "ravine" functions and the Minkowski-Cohn surface}\\

The closed Minkowski-Cohn surface is defined by a corresponding "ravine" function and is non-differentiable. However, it is possible to construct corresponding functions that are continuous on open subdomains of the domain of definition of the Minkowski-Cone surface. Let us define such functions—evaluated on this surface—as well as the corresponding derivatives of the Minkowski-Cone surface (including the functions themselves). For instance, the second derivative with respect to \(\sigma\), evaluated at \(\sigma_{p}\), is denoted as \(\Delta_{\sigma^{2}}^{''}|_{\sigma = \sigma_{p}}\).

We have: \(\Delta(p, 1) = \Delta_{p}^{(1)} = 4^{- \frac{1}{p}}\frac{1 + \tau_{p}}{1 - \tau_{p}}\), \(\Delta(p, \sigma_{p}) = \Delta_{p}^{(0)} = \frac{1}{2}\sigma_{p}\).\\

For certain issues, the following known asymptotic expansions may be useful.

Expansions of the functions \(\Delta_{p}^{(0)}\) and \(\Delta_{p}^{(1)}\) in the neighborhoods of \(p = 1\), \(p = 2\), and \(p = \infty\):\\

Denote by $\Delta ( D_p)$ the critical determinant of the domain $ D_p$.

\begin{equation}
\label{cdd}
\Delta(D_p) = \left\{
                   \begin{array}{lc}
     {\Delta^{(0)}_p}( D_p), \; 1 \le p \le 2, \; p \ge p_{0},\\
    {\Delta^{(1)}_p}( D_p), \;  2 \le p \le p_{0};\\
                     \end{array}
                       \right.
 \end{equation}

\textbf{Lemma 1}. \(\sigma_{p}=1+2(p-1)\log 2-(p-1)^{2}\bigl(2\log 2-\log^{2}2\bigr)+O((p-1)^{3})\).\\

\textbf{Corollary 1}.\\

\[\Delta_{p}^{(0)} = \frac{\sigma_{p}}{2} \approx \frac{1}{2}\left\{ 1 + 1.3862(p - 1) - 0.9059(p - 1)^{2} \right\}.\]\\

\textbf{Lemma 2.} Given that \(p\) is fixed, and \(\sigma\) increases from \(1\) to \(\sigma_{p} = (2^{p} - 1)^{1/p}\):

1) \((1 + \sigma^{p})^{- 1/p}\) decreases from \(2^{- 1/p}\) to \(\frac{1}{2}\).

2) \(\sigma(1 + \sigma^{p})^{- 1/p}\) increases from \(2^{- 1/p}\) to \(\frac{1}{2}\sigma_{p}\);

3) \(\tau\) decreases from \(\tau_{p}\) to 0;

4) \((1 + \tau^{p})^{- 1/p}\) increases from \((1 + \tau_{p}^{p})^{- 1/p}\) to 1;

5) \(\tau(1 + \tau^{p})^{- 1/p}\) decreases from \(\tau_{p}(1 + \tau_{p}^{p})^{- 1/p}\) to \(0\); 6) \emph{A} increases from \(\tau_{p}(1 + \tau_{p}^{p})^{- 1/p}\) to \(1/2\);

7) \emph{B} decreases from \((1 +  \tau_{p}^{p})^{- 1/p}\) to \(\frac{1}{2}\sigma_{p}\).\\

\textbf{Lemma 3.} For any \(p\) and \(\sigma\)

\(\tau < 1 - (1 +  \sigma^{- p})^{- 1/p}\).

\textbf{Corollary 2.} For any \(p\) and \(\sigma\)

\(\tau < \frac{\log(1 + \sigma^{-p})}{p} < \frac{1}{p\sigma^{p}}\).

\textbf{Corollary 3.} For \(p > 1\)

\(\tau \leq \tau_{p} < 1 - 2^{- 1/p} < \frac{\log 2}{p} < \frac{0.7}{p}\).\\

\textbf{1.1.8 Iterative calculation} of \(\tau_{p}\) \textbf{as an implicit function defined by the equation} \(2\left( 1 - \tau_{p} \right)^{p} = 1 + \tau_{p}^{p},\ 0 < \tau_{p} < 1.\)\\

\textbf{Lemma 4.} The function \emph{F(}\(\tau_{p}\)\emph{) =} \(2\left( 1 - \tau_{p} \right)^{p} - 1 -  \tau_{p}^{p}\) is monotonically decreasing with respect to \(\tau_{p}\) on the half-open interval (\emph{0.1; 0.7{]}} and changes sign within this interval.\\

\textbf{Proof.} Indeed, \(F_{\tau_{p}}^{'} = - 2p(1 - \tau_{p})^{p - 1} - p\tau_{p}^{p - 1} < 0\), which follows from point 3) of Lemma 2, thereby establishing monotonicity. the sign change of the \emph{function F(}\(\tau_{p}\)\emph{)} on the half-interval (\emph{0.1; 0.7{]}} is obvious.\\

\textbf{Procedure} \(\tau_{p}\) \textbf{calculation of dichotomous division:}\\

\begingroup
\parindent=0pt\parskip=0pt\relax
\fontsize{10.3pt}{12.8pt}\selectfont
\raggedright\ttfamily\sloppy
\noindent\texttt{TPD: \texttt{Proc} (p,}\(\tau_{p}\)\texttt{);}\par
\noindent\texttt{x0 = 0; x1 = 0.7/p; E1 = 10-9;}\par
\noindent\texttt{M: t1 = x0;}\par
\noindent\hspace*{18.2pt}\texttt{F0 =2(1 - t1}\()^{p}\) \texttt{- 1 - t1p;}\par
\noindent\hspace*{34.0pt}\texttt{t1 = x1;}\par
\noindent\hspace*{18.2pt}\texttt{F1 =2(1 - t1}\()^{p}\) \texttt{- 1 - t1p;}\par
\noindent\hspace*{34.0pt}\texttt{x2 = (x0 + x1)/2; t1 = x2;}\par
\noindent\hspace*{18.2pt}\texttt{F2 =2(1 - t1}\()^{p}\) \texttt{- 1 - t1p;}\par
\noindent\hspace*{11.3pt}\texttt{\texttt{if} ABS(F2) \textless{} E1, \texttt{then} \texttt{goto} M1;}\par
\noindent\hspace*{25.0pt}\texttt{\texttt{if} F0F1 \textless{} 0, \texttt{then} \texttt{do}; x0 = x2; \texttt{goto} M; \texttt{end};}\par
\noindent\texttt{M1:} \(\tau_{p}\) \texttt{= x2;}\par
\noindent\texttt{\texttt{end} TPD;}\par
\endgroup

We present selected calculation results for this program, summarized in the table below. In the table, a minus sign (--) following a number indicates that the value is an underestimate (truncated/rounded down), while a plus sign (+) indicates an overestimate (rounded up).\\

\begin{flushright}Table 1\end{flushright}


\begin{center}
{\setlength{\tabcolsep}{18pt}\renewcommand{\arraystretch}{1.14}
\begin{tabular}{|c|c|}
\hline
$p$ & $\tau_p$ \\ \hline
1.01 & 0.332751- \\ \hline
1.1 & 0.327427+ \\ \hline
1.5 & 0.302104+ \\ \hline
2 & 0.267949- \\ \hline
2.3 & 0.247392+ \\ \hline
2.5 & 0.234162+ \\ \hline
\end{tabular}
}
\end{center}




\textbf{1.1.8. Iterative computation} of \(\tau = \tau(p,\sigma)\) \textbf{as an implicit function defined by} \(A^{p} + B^{p} = 1, 0 \leq \tau \leq \tau_{p}\). The equation \(A^{p} + B^{p} = 1\), along with its variables and their respective domains, is described in Section 2. The value \(\tau = \tau(p,\sigma)\) is computed either at a single point \((p, \sigma)\) within the aforementioned domain, or across a finite set of points within this domain—typically along straight lines. For the computation of \(\tau = \tau(p,\sigma)\) presented below, we will introduce the bisection method, the method of chords, and a technique we shall refer to as the Malyshev iterative method.\\




\textbf{Lemma 4.} Given that \(p\) is fixed and $\sigma$ increases from \(1\) to \(\sigma_{p} = (2^{p} - 1)^{1/p}\), the function $\tau$ monotonically decreases from \(\tau_{p}\) to 0.\\

\textbf{Proof.} Under the conditions of the lemma,

\[\frac{\partial\tau}{\partial\sigma} = - \frac{(1 + \tau^{p})^{1 + 1/p}(B^{p - 1} + A^{p - 1}\sigma^{p - 1})}{(1 + \sigma^{p})^{1 + 1/p}(B^{p - 1} - A^{p - 1}\tau^{p - 1})}.\]

It follows from Lemma 1 that both the numerator and the denominator of the fraction (the term following the minus sign) are greater than zero. This implies that $\tau$ is monotonically decreasing with respect to the variable $\sigma$.\\

\textbf{1.2 Optimal packing of unit Minkowski spheres in the plane}\\

A well-known monograph \cite{src15} presents and investigates a wide range of problems and methods concerning the theory of lattices, groups, and sphere packing. Recently, M. Viazovska and her colleagues obtained remarkable results regarding the optimal packing of spheres in dimensions 8 and 24 \cite{src16,src17}. The cited works investigate the packing of Euclidean spheres.\\

In his monograph \cite{src5}—as well as in earlier works (see references in \cite{src5})—Minkowski, while studying problems in the theory of Diophantine approximations, defined objects that we refer to in this article as Minkowski spheres.\\

Let \(|a|\) denote the absolute value of a real number \(a\). The \(n\)-dimensional unit Minkowski sphere corresponding to a parameter \(p \geq 1\), centered at the point \((a_{1},\ldots,a_{n})\), is defined as follows:\\

\[S_{p}^{n - 1}(a_{1},\ldots,a_{n}):\left| x_{1} - a_{1} \right|^{p} + \cdots + \left| x_{n} - a_{n} \right|^{p} = 1.\]\\

Minkowski himself defined these objects for the case \(n = 2\).\\

The 1-dimensional unit Minkowski sphere corresponding to a parameter \(p \geq 1\) (and we consider only such spheres in the article \cite{src24}), centered at the point \((a_{1},a_{2})\), is defined as follows:

\[S_{p}^{1}(a_{1},a_{2}):\left| x_{1} - a_{1} \right|^{p} + \left| x_{2} - a_{2} \right|^{p} = 1.\]\\

\textbf{1.2.1. Interior and Exterior Points of a Body}\\

Let us recall the standard definition of interior and exterior points of a geometric body, using the one-dimensional Minkowski spheres we are studying as an example.\\

\textbf{Definition.} \emph{Consider a one-dimensional Minkowski sphere in the real plane, defined by the equation}\\

\[S_{p}^{1}(a_{1},a_{2}):\left| x_{1} - a_{1} \right|^{p} + \left| x_{2} - a_{2} \right|^{p} = 1.\]\\

\emph{centered at a point \(a = (a_{1}, a_{2})\) in this plane. Let there be a point \(b = (b_{1}, b_{2})\), and let a line \(L\) intersect the Minkowski sphere at points \(c = (c_{1}, c_{2})\) and \(d = (d_{1}, d_{2})\). By the construction of the Minkowski sphere, the line segments connecting point \(a\) to point \(c\) and to point \(d\) are equal in length.} Let us denote these segments as $[a,c]$ and $[a,d]$, respectively, where $[a,c] = [a,d]$.

A Minkowski sphere positioned in this way divides the points of the plane into: *interior points*—those whose distance from the center of the sphere is less than the length of the segment $[a,c]$; *exterior points*—those whose distance from the center is greater than the length of the segment $[a,c]$; and points *on* the Minkowski sphere itself—those whose distance from the center is equal to the length of the segment $[a,c]$.\\

\textbf{Remark.} Interior and exterior points do not belong to the Minkowski sphere.\\

\textbf{1.2.2. On the Packing of Minkowski Spheres}\\

If the center of a Minkowski sphere coincides with the origin, we use the notation $S_{p}^{1}$ and refer to such a sphere as centrally symmetric.\\

Note that there is (up to rigid motions in the Euclidean plane) only one unit centrally symmetric Euclidean sphere, whereas there exists a continuum of unit Minkowski spheres, parameterized by a parameter $p \geq 1$. In what follows, we consider only unit Minkowski spheres.

We say that a system of identical Minkowski spheres in a 2-dimensional space forms a packing if no two spheres in the system share any interior points.\\

Packing is related to packing density and is considered optimal when the density is maximized. Recall that the optimal packing of bodies in a plane is a lattice packing \cite{src22,src23}. The optimal packing of Minkowski unit spheres relies on results derived from the proof of the Minkowski conjecture \cite{src9}.\\

In a previous note \cite{src10}, I provided a proof of the theorem regarding the optimal packing of Minkowski unit balls, formulated in terms of Minkowski balls. The relationship between Minkowski spheres and Minkowski balls is straightforward: a Minkowski sphere is the boundary of a Minkowski ball. In this paper, I present a detailed proof—conceptually simple—formulated in terms of Minkowski spheres (which are more relevant to applications of packing, such as in coding theory). I also define a partition of Minkowski spheres into three classes, the names of which reflect the history of research into the Minkowski conjecture \cite{src7,src8,src9,src10}. These three classes—whose union constitutes a partition of the set of Minkowski spheres and which emerge from the proof of the Minkowski conjecture—are designated, respectively, as Watson spheres, Davis spheres, and Mordell-Chebyshev spheres.\\

Our main result is a theorem concerning optimal lattice packings in the real plane for Watson spheres, Davis spheres, and Mordell-Chebyshev spheres. We prove that such optimal packings are determined by sublattices of index two of the corresponding critical lattices.\\

Such lattices can be computed by applying the method of parameterizing admissible lattices—within the framework of Minkowski’s conjectures—to doubled Minkowski spheres and their admissible lattices.\\

\textbf{1.2.2 Terminology and the Result Regarding the Proof of the Minkowski Conjecture.}\\

We shall establish our terminology based on \cite{src1,src11,src12}. Definitions are provided for an \(n\)-dimensional real space, with specific refinements for the case of the real plane—the scenario under investigation here. We define a \emph{body} (or \emph{region}) as a connected set \(M\) that possesses interior points and is contained within the closure of its open kernel.

Let \(\mathcal{R}\) be a body, and let \(\Lambda\) be a full lattice in \({\mathbb{R}}^{n}\) with a basis \(\left\{ a_{1},\ldots,a_{n} \right\}\). The lattice \(\Lambda\) is said to be \emph{admissible} for the body \(\mathcal{R}\) (or \(\mathcal{R}\)-admissible) if the intersection of \(\Lambda\) and \(\mathcal{R}\) is either empty or consists solely of the origin. Let \(\overline{\mathcal{R}}\) denote the closure of \(\mathcal{R}\). A lattice \(\Lambda\) is strictly admissible for \(\overline{\mathcal{R}}\) (\(\overline{\mathcal{R}}\) being strictly admissible) if the intersection of \(\overline{\mathcal{R}}\) and \(\Lambda\) is either empty or consists solely of the origin.\\

Let \(d(\Lambda) = \left| \det\left( a_{1},\ldots,a_{n} \right) \right|\) be the determinant of the lattice \(\Lambda\).\\

The infimum (greatest lower bound) \(\Delta(\mathcal{R})\) of the determinants of all lattices admissible for \(\mathcal{R}\) is called the critical determinant of the body \(\mathcal{R}\); if there are no lattices admissible for \(\mathcal{R}\), we set \(\Delta(\mathcal{R}) = \infty\).\\

A lattice \(\Lambda\) is called critical (for a body \(\mathcal{R}\)) if \(d(\Lambda) = \Delta(\mathcal{R})\).\\

For the Minkowski sphere \(S_{p}^{1}\), the critical determinant of the interior \(D_{p}\) of the sphere \(S_{p}^{1}\) is referred to as its critical determinant and is denoted by \(\Delta(D_{p})\) \cite{src9}.\\

Using the notation from \cite{src9}, the following result was proven.\\

\textbf{Theorem} 1 \cite{src9}.

\begin{equation}\label{eq:critical-determinant}
\Delta(D_p)=
\begin{cases}
\Delta(p,1), & 1\leq p\leq 2,\quad p\geq p_0,\\[2pt]
\Delta(p,\sigma_p), & 2\leq p\leq p_0;
\end{cases}
\end{equation}
\[\Delta\left( p_{0},\sigma_{p} \right) = \Delta\left( p_{0},1 \right),\ 2.57 < p_{0} < 2.58,\ p_{0} \approx 2.5725.\]

The constant \(p_{0}\) is called the Davis constant.\\

\textbf{Corollary.} Using the notation from \cite{src9}, we have the following expressions for the corresponding critical determinants and lattices:

\[\Delta_{p}^{(0)} = \Delta\left( p,\sigma_{p} \right) = \frac{1}{2}\sigma_{p},\ \sigma_{p} = \left( 2^{p} - 1 \right)^{1/p},\]

\[\Delta_{p}^{(1)} = \Delta(p,1) = 4^{- \frac{1}{p}}\frac{1 + \tau_{p}}{1 - \tau_{p}}, \quad 2{(1 - \tau_{p})}^{p} = 1 + \tau_{p}^{p},\ 0 \leq \tau_{p} < 1.\]\\

For these critical determinants, the corresponding critical lattices possess the following properties: they are two \(D_{p}\)-admissible lattices, and each contains three pairs of points on \(S_{p}^{1}\) such that

\[(1,0) \in \Lambda_{p}^{(0)},\quad \left( {- 2}^{- \frac{1}{p}},2^{- \frac{1}{p}} \right) \in \Lambda_{p}^{(1)}.\]\\

Moreover, these conditions uniquely determine the corresponding lattices.\\

\textbf{1.2.3. Minkowski Sphere Packings and Partitions of the Set of One-Dimensional Minkowski Spheres}\\

Recall the definition of a packing lattice. We present it here for \emph{n}-dimensional Minkowski spheres \(S_{p}^{n}\).\\

\textbf{Definition 1}. Let \(\Lambda\) be a full-rank lattice in \({\mathbb{R}}^{n}\), and let \(a\) be an \(n\)-dimensional vector in \({\mathbb{R}}^{n}\). If no two spheres in the collection \(\{S_{p}^{n} + b \mid b \in \Lambda + a\}\) share any interior points, then the collection \(\{S_{p}^{n} + b \mid b \in \Lambda + a\}\) is called an \((S_{p}^{n}, \Lambda)\)-packing, and \(\Lambda\) is called a packing lattice for \(S_{p}^{n}\).\\

Recall also that if \(\alpha \in {\mathbb{R}}\) and \(D_{p}^{n}\) is a set in \({\mathbb{R}}^{n}\), then \(\alpha D_{p}^{n}\) denotes the set of points \(\{\alpha x \mid x \in D_{p}^{n}\}\).\\

\textbf{Proposition 1}. A lattice \(\Lambda\) is a packing lattice for \(S_{p}^{1}\) if and only if it is an admissible lattice for \(2 \times S_{p}^{1}\).\\

\textbf{Proof}. A lattice admissible for \(S_{p}^{1}\) is strictly admissible for the interior \(D_{p}\) of the sphere, while a lattice admissible for \(2 \times S_{p}^{1}\) is strictly admissible for \(2 \times D_{p}\). Let \(\Lambda\) be a packing lattice for \(S_{p}^{1}\) that is not a strictly admissible lattice for \(2 \times D_{p}\). Then \(\Lambda\) contains a point \(a\) distinct from the origin. It follows that \(D_{p}\) and \(D_{p} + a\) share a common point, \(\frac{1}{2}a\). Consequently, \(\Lambda\) is not a packing lattice for \(S_{p}^{1}\).

Now suppose that \(\Lambda\) is not a packing lattice for \(D_{p}\). Then there exist two distinct points \(b_{1}, b_{2} \in \Lambda\) and a point \(c\) such that \(c \in D_{p} + b_{1}\) and \(c \in D_{p} + b_{2}\). Thus, there exist points \(a_{1}, a_{2} \in D_{p}\) such that \(c = a_{1} + b_{1} = a_{2} + b_{2}\). Therefore, \(b_{1} - b_{2} = a_{2} - a_{1} \in D_{p}\), while \(b_{1} - b_{2} \neq 0\) and \(b_{1} - b_{2} \in \Lambda\).

Therefore, \(\Lambda\) is not a (strictly) admissible lattice.

\textbf{Proposition 2}. The density of a \((D_{p}, \Lambda)\)-packing is equal to \(V(D_{p})/d(\Lambda)\) and is maximal if \(\Lambda\) is critical for \(2 \times D_{p}\).\\

Proposition 2 follows from general results \cite{src15,src22,src23} regarding the density of packings of convex centrally symmetric bodies and their closures, applied in our case to 2-dimensional Minkowski spheres.\\


 We have:\\

Watson spheres: \(|x|^{p} + |y|^{p} = 1,  1 < p < 2;\)\\

Davis spheres: \(|x|^{p} + |y|^{p} = 1,  2 \leq p < p_{0};\)\\

Mordell-Chebyshev spheres: \(|x|^{p} + |y|^{p} = 1,  p \geq p_{0};\)\\

limiting Watson sphere \(S_{1}^{1}:|x| + |y| = 1;\)\\

limiting Mordell-Chebyshev sphere \(S_{\infty}^{1}:\left\|x,y\right\| = \max(|x|,|y|)\).\\

\textbf{1.2.4. Optimal Packing of Watson Spheres, Davis Spheres, and Mordell-Chebyshev Spheres in the Plane}.\\

We consider possible optimal lattice packings of these spheres and their connection to critical lattices.\\

First, we present the lattices corresponding to trivial optimal lattice packings for the limiting (asymptotic) cases at \(p = 1\) and \(p = \infty\) ("infinity"), followed by an introductory example: the optimal packing of one-dimensional unit spheres.\\

\textbf{Examples}.

1) The lattice

\[\Lambda_{1}^{(1)} = \left\{ \left( \frac{1}{2},\frac{1}{2} \right),(0,1) \right\}\]

is a critical lattice for \(S_{1}^{1}\). The limiting case of Minkowski spheres for \(p = 1\) (limiting Watson spheres) yields an optimal sublattice of index two within the lattice \(\Lambda\)—a packing lattice with density 1. In this case, the centers of the Watson spheres are located at the vertices of a sublattice of index two within the lattice \(\Lambda\).

2) The lattice

\[\Lambda_{\infty}^{(1)} = \left\{ (1,1),(0,1) \right\}\]

is a critical lattice for \(S_{\infty}^{1}\). The limiting case of Minkowski spheres for $p = \infty$ (Mordell-Chebyshev spheres) yields an optimal sublattice of index two within the lattice $\Lambda$—a packing lattice with density 1. In this case, the centers of the Mordell-Chebyshev spheres are located at the vertices of a sublattice of index two within the lattice $\Lambda$.

3) The lattice

\[\Lambda_{2}^{(0)} = \left\{ (1,0),\left( \frac{1}{2},\frac{\sqrt{3}}{2} \right) \right\}\]

is a critical lattice for $S_{2}^{1}$. The Davis lattice sphere packing for $p = 2$ yields an optimal packing with a density of $\approx 0.91$ relative to a sublattice of index two within the critical lattice $\Lambda$. In this case, the centers of the Davis spheres are located at the vertices of a sublattice of index two within the lattice $\Lambda$.

\textbf{Proposition 3}. If $\Lambda$ is a critical lattice for the Minkowski sphere $S_{p}^{1}$, then a sublattice $\Lambda_2$ of index two within this critical lattice is a critical lattice for $2 \times S_{p}^{1}$.

(See the examples for $p = 1, 2, \infty$ above).

\textbf{Proof}. Since the Minkowski sphere $S_{p}^{1}$ is symmetric with respect to the origin and is convex, $2 \times S_{p}^{1}$ is also convex and symmetric with respect to the origin \cite{src23}. In the parameterization, admissible lattices \(\Lambda\)—which possess three pairs of points on the boundary of the sphere \(S_{p}^{1}\)—are parameterized by points on a manifold \cite{src5} as follows:

\begin{equation}\label{eq:admissible-lattice}
 \Lambda = \left\{ \left( \left( 1 + \tau^{p} \right)^{- \frac{1}{p}},\tau\left( 1 + \tau^{p} \right)^{- \frac{1}{p}} \right),\left( - \left( 1 + \sigma^{p} \right)^{- \frac{1}{p}},\sigma\left( 1 + \sigma^{p} \right)^{- \frac{1}{p}} \right) \right\}
\end{equation}
where \(0 \leq \tau < \sigma\) and \(0 \leq \tau \leq \tau_{p}\), with

\[{2\left( {1 - \tau}_{p} \right)}^{p} = 1 + \tau_{p}^{p}, \quad 0 \leq \tau_{p} < 1,\]

\[1 \leq \sigma \leq \sigma_{p}, \quad \sigma_{p} = \left( 2^{p} - 1 \right)^{\frac{1}{p}}.\]

Admissible lattices for doubled Minkowski spheres \(2 \times S_{p}^{1}\) take the form

\begin{equation}\label{eq:doubled-lattice}
 \Lambda_{{2D}_{p}} = \left\{ 2\left( \left( 1 + \tau^{p} \right)^{- \frac{1}{p}},2\tau\left( 1 + \tau^{p} \right)^{- \frac{1}{p}} \right),\left( - 2\left( 1 + \sigma^{p} \right)^{- \frac{1}{p}},2\sigma\left( 1 + \sigma^{p} \right)^{- \frac{1}{p}} \right) \right\}
\end{equation}
Thus, the moduli space for these admissible lattices is given by

\begin{equation}\label{eq:doubled-determinant}
\Delta(p,\sigma)_{{2D}_{p}} = 4(\tau + \sigma)\left( 1 + \tau^{p} \right)^{- \frac{1}{p}}\left( 1 + \sigma^{p} \right)^{- \frac{1}{p}}
\end{equation}
over the same domain

\(\ M:\infty > p > 1,\ 1 \leq \sigma \leq \sigma_{p} = \left( 2^{p} - 1 \right)^{- \frac{1}{p}}\)

as the moduli space for \eqref{eq:admissible-lattice}. We refer to the moduli space \eqref{eq:doubled-determinant} as the Minkowski–Cohn manifold of the admissible lattices \eqref{eq:doubled-lattice}. Thus, the critical determinants of doubled Minkowski spheres have the following representations:

\[\Delta_{p}^{(0)}\left( 2D_{p} \right) = \Delta\left( p,\sigma_{p} \right)_{{2D}_{p}} = 2  \sigma_{p},\]\\

\[\sigma_{p} = \left( 2^{p} - 1 \right)^{- \frac{1}{p}},\]\\

\[\Delta_{p}^{(1)}\left( 2D_{p} \right) = \Delta(p,1)_{{2D}_{p}} = 4^{1 - \frac{1}{p}}\frac{1 + \tau_{p}}{1 - \tau_{p}},\]\\

\[2\left( 1 - \tau_{p} \right)^{p} = 1 + \tau_{p}^{p},\]  \[0 \leq \tau_{p} < 1.\]\\

\textbf{Theorem 2.} The optimal lattice packing of Watson spheres, Davis spheres, and Mordell-Chebyshev spheres is realized with respect to sublattices of index two of the critical lattices

\[(1,0) \in \Lambda_{p}^{(0)},\left( {- 2}^{- \frac{1}{p}},2^{- \frac{1}{p}} \right) \in \Lambda_{p}^{(1)}.\]

\textbf{Proof.} By Proposition 3, the critical lattice \(2 \times S_{p}^{1}\) is a sublattice of index two of the critical lattice of the Minkowski sphere \(S_{p}^{1}\). Thus, it is an admissible lattice for \(2 \times S_{p}^{1}\), and according to Proposition 1 is a packing lattice \(S_{p}^{1}\). According to Proposition 2, the corresponding lattice packing has maximum density and is, therefore, optimal.

\newpage

\textbf{References}



\end{document}